\documentclass{article}
\usepackage{amsmath,amssymb,amsthm}
\newtheorem{thm}{Theorem}
\newtheorem{lem}[thm]{Lemma}
\newtheorem{prop}[thm]{Proposition}
\theoremstyle{definition}
\newtheorem{dfn}[thm]{Definition}

\def\beq #1 \e{\begin{equation}#1\end{equation}}
\def\bml #1 \e{\begin{multline}#1\end{multline}}
\def\bmlg #1 \e{\begin{multline*}#1\end{multline*}}
\let\SS\S
\def\cf{coefficient}
\def\f{function}
\def\ip{interpolation}
\def\r{remainder}
\def\bL{\mathbb L}
\def\bW{\mathbb W}
\def\cH{\mathcal H}
\def\cL{\mathcal L}
\def\cM{\mathcal M}
\def\cR{\mathcal R}
\let\D\Delta
\def\R{\mathbb R}
\let\S\Sigma
\let\er\eqref
\let\d\delta
\let\ve\varepsilon
\let\vf\varphi
\def\fint#1{\mathchoice
{\hbox{\ooalign{\hfil$\relbar$\hfil\crcr$\displaystyle\intop_{#1}$}}}
{\mathop{\hbox{\ooalign{\hfil$\relbar$\hfil\crcr$\textstyle\intop$}\!}}\nolimits_{#1}}
{\mathop{\hbox{\ooalign{\hfil$\relbar$\hfil\crcr$\textstyle\intop$}\!}}\nolimits_{#1}}
{\mathop{\hbox{\ooalign{\hfil$\relbar$\hfil\crcr$\textstyle\intop$}\!}}\nolimits_{#1}}}
\def\loc{_{\rm loc}}
\let\n\nabla
\let\pa\partial
\let\q\quad
\let\s\sigma
\let\ti\textit
\let\wh\widehat
\let\wt\widetilde

\numberwithin{equation}{section}

\begin{document}
\author{B. Bojarski\thanks{Partially supported by the Polish Ministry of Science
grant no.\ N N201 397837 (years 2009--2012) and the Academy of Finland.}\\
Institute of Mathematics, Polish Academy of Sciences\\
00-956 Warszawa, Poland\\
E-mail: b.bojarski@impan.pl}
\title{Sobolev spaces and Lagrange interpolation}
%\date{Preliminary version, \today}
\maketitle

In this short paper the discussion of the pointwise characterization of
functions $f$ in the Sobolev space $W^{m,p}(\R^n)$ given in the recent paper
\cite{r0} is supplemented in \SS1 by a direct, essentially geometric, proof of
the novel inequality (for $m>1$), appearing in \cite{r0} apparently for the
first time, and involving the use of the $m$-th difference of the
function~$f$. 

Moreover in \SS2 some additional comments to the text in \cite{r0} are given
and a natural class of Sobolev spaces in domains $G$ in $\R^n$ is defined. 
\SS3~contains some final remarks.

\section{}
Let us recall \cite{r8}, \cite{r7} that for an arbitrary integer $l\ge0$ and a
real or complex valued \f\ $f$ on $\R^n$ the expression
\beq\label{e1.1}
\D_h^l f(x):=\sum_{j=0}^l (-1)^{l-j}\binom lj f(x+jh)=
(-1)^l \sum_{j=0}^l (-1)^j \binom lj f(x+jh)
\e
is called the $l$-th difference of the \f\ $f$ at the point $x\in\R^n$ with
step~$h$, $h\in\R^n$, $h\ne0$. We set also $\D_h^0 f(x):=f(x)$ and
$\D_0^l f(x):=0$. As is classically known, for $y=x+lh$, $h=\frac{y-x}l$,
\er{e1.1} has a beautiful interpretation as the difference
or error in approximating the \f\ $f(y)$ by its interpolating polynomial
evaluated at~$y$
\bml\label{e1.2}
\D^lf(x;y)=\D_h^l f(x)=f(y)-\cL(y;f;x_0,\dots,x_{l-1}) 
\equiv (-1)^l\wt\D{}^l f(x,y)\\ \quad (x_0=x),
\e
where $\cL(y;f;x_0,\dots,x_{l-1})\equiv \sum_{j=0}^{l-1} f(x_j)\ell_j
(y,x_0,\dots,x_{l-1})$ is the Lagrange interpolating polynomial for the
\f~$f$ and the equidistant colinear nodes $x_i=x_0+ih$, $i=0,\dots,l-1$,
\cite{r7}. 
Here $\ell_j(y,x_0,\dots,x_{l-1})$ stand for the fundamental Lagrange
polynomials in~$y$. $\wt\D{}^l(x,y)$ is the notation used in \cite{r0}. Let us remark
that all points $x_j=x+jh$, $j=0,\dots,l-1$, $x_l=y$, are situated on the
affine line $\cR$ in~$\R^n$, joining $x$ and~$y$, which can be identified
with the real line $\R^1$, making all the algebraic operations inherent in
\er{e1.2} meaningful.

In the sequel we shall use for \er{e1.2} the term \ti{Lagrange interpolation
remainder}, or just \ti{Lagrange remainder}, of order~$l$ at the point $x$ 
evaluated at~$y$, in
analogy with the term \ti{Taylor--Whitney remainder} centered at~$x$,
\beq\label{e1.3}
R^{l-1}f(x;y):=f(y)-T_x^{l-1}f(y),
\e
now common it mathematical literature.
In \cite{r0} $\wt\D{}^l f(x,y)$ was also called an $l$-th finite difference remainder
of the \f~$f$ at~$x$ evaluated at~$y$.

For \f s $f\in W^{m,p}(\R^n)$, $p>1$, the fundamental novel inequality
referred to above reads as
\beq \label{e1.4}
|\D^m f(x;y)| \le |x-y|^m [\wh a_f^m(x)+\wh a_f^m(y)]
\e
for some $\wh a_f\in L^p(\R^n)$.

We skip here over the somewhat delicate point that the Sobolev \f s in
general do not have pointwise values and the left hand side of inequality
\er{e1.4} is meaningful only up to subsets of measure zero. The right hand
side may be infinite on a non-empty set of measure zero.

The \f al \cf s $\wh a_f^m(x)$ in \er{e1.4} are not uniquely defined. They
are collectively called \ti{mean maximal $m$-gradients} of the \f~$f$ and
play the role of a variable Lipschitz \cf\ of~$f$. Roughly speaking, they all
can be majorized by the local maximal \f\ of the generalized Sobolev gradient
$|\n^m f| $ of $f \in W^{m,p}(\R^n)$ as will be also seen from the constructive
proof of~\er{e1.4} sketched below. 

Our proof is organized in a series of lemmata.

\begin{lem}\label{l1}
For $f\in W^{1,p}(\R^n)$, $1<p<\infty$, the following inequality holds
\beq \label{e1}
|f(x)-f(y)|\le |x-y|\bigl(a_f^\d(x)+a_f^\d(y)\bigr), \q x,y\in\R^n
\e
for some $a_f^\d\in L\loc^p(\R^n)$, $\d=|x-y|$.
\end{lem}

\proof
For arbitrary $x,y\in\R^n$ we have
\beq \label{e2}
f(x)-f(y)=\int_0^1\langle \n f(x+ht),h\rangle\,dt, \q h=y-x,
\e
hence
\beq \label{e2'}
|f(x)-f(y)|\le |x-y|\int_0^1|\n f|(x+ht)\,dt.
\e

Let $B(x,r)$ be the ball of radius $r$ centered at $x$, and $\S_r(x,y)$ the
spherical segment
\[
\S_r=B(x,r)\cap B(y,r), \q r=|x-y|.
\]
For an arbitrary $z\in\S_r$
\beq \label{e3}
|f(x)-f(y)|\le |f(x)-f(z)|+|f(z)-f(y)|.
\e
Since $|x-z|\le|x-y|$, $|y-z|\le |x-y|$, averaging \er{e3} over $z\in\S_r$ 
we get
\bml
|f(x)-f(y)|\le \fint{\S_r}|f(z)-f(x)|\,d\s_z
+\fint{\S_r}|f(z)-f(y)|\,d\s_z\\
{}\le \frac{|B(x,r)|}{|\S_r|}\Bigl(\fint{B(x,r)}|f(z)-f(x)|\,d\s_z
+\fint{B(y,r)}|f(z)-f(y)|\,d\s_z\Bigr),
\e
where the notation $|G|$ for a subset $G$ in $\R^n$ is used for the volume of
$G$, $|B(x,r)|=|B(y,r)|$.

By elementary geometry the ratio $\frac{|B(x,r)|}{|\S_r|}$ is a constant
depending only on~$n$, $\frac{|B(x,r)|}{|\S_r|}=C(n)$, and, as is well known,
the average $\fint{B(x,r)}|f(z)-f(x)|\,d\s_z$ is estimated by the local
Hardy-Littlewood maximal \f\ at~$x$ of the gradient $|\n f|$, $\cM^\d(|\n
f|)(x)$, \cite{St1}, \cite{St2}. 
Here the inequality \er{e2'} is used. Thus in~\er{e1} the \f\
$a_f^\d(x)$ is controlled by $\cM^\d(|\n f|)(x)$: in fact $a_f^\d(x)
\le C(n)\cM^\d(|\n f|)(x)$.
\endproof

The proof above, without changes, works for vector valued \f s. This proof
should be compared with the proof of the basic pointwise inequality (1) in
our paper with P.~Haj\l asz from 1993 \cite{r1}. Notice that it does not
refer to Riesz potentials and Hedberg lemma as in~\cite{r1}. 

When combined
with Reshetnyak's trick \cite{r4} used in~\cite{r1}, it can be used to deduce,
in a direct way, a new and simplified proof of the basic pointwise
inequalities in \cite{r3, r1, r2}. 

Let $f\in W^{k,p}(\R^n)$. For $l=0,1,\dots,k-1$, $h\in\R^n$, consider the \f
s 
\beq \label{e4}
g_h^l(x)=\underbrace{\int_0^1 \cdots \int_0^1}_{l\ \rm times}
\n^l f\Bigl(x+\sum_{i=1}^l t_i h\Bigr)(h,\dots,h)\,dt_1 \dots dt_l,
\e
where $\n ^lf$ is the $l$-th gradient of~$f$ considered as an $l$-polylinear
form on~$\R^n$.

\begin{lem} \label{l2}
$g_h^l(x)$ as a \f\ of $x\in\R^n$ is in the class $W^{k-l,p}(\R^n)$.
\end{lem}

\proof Obvious. \endproof

In particular, $g_h^0(x)\equiv f(x)$,
\[
g_h^1(x)=\int_0^1\langle\n f(x+th,h\rangle \,dt
\equiv \int_0^1 \sum_{i=1}^n \frac{\pa f}{\pa x_i}(x+ht)h_i\,dt, \q \hbox{etc.}
\]

\begin{lem}\label{l3}
The \f\ $g_h^l(x)$ for $l=0,\dots,k-1$ has the integral representation
\bml \label{e5}
g_h^l(x)=\sum_{j=0}^l (-1)^l \binom lj f(x+jh)\\ {}=
\int_0^1 \cdots \int_0^1 \n^l f\Bigl(x+\sum_{i=1}^l t_i h\Bigr)
(\underbrace{h,\dots,h}_{l\ \rm times})\,dt_1 \dots dt_l.
\e
\end{lem}

\proof
This is the well known formula of finite difference calculus \cite{r8, r7}.
For $l=k$ it is also used in the paper \cite{r5} of R. Borghol.
\endproof

\begin{lem}\label{l4}
For $l=k-1$ we have the formula
\bml \label{e6}
g_h^{k-1}(x)-g_h^{k-1}(x+h)=
\sum_{l=1}^k (-1)^{l-1}\binom{k-1}{l-1}\bigl[f(x+(l-1)h)-f(x+lh)\bigr]\\
{}\equiv \sum_{j=0}^k (-1)^j \binom kj f(x+jh)=\wt \D{}^k f(x,y)
=(-1)^k\D^k f \q \hbox{for }h=\frac{y-x}{k},
\e
in the notation of {\rm\cite{r0}}.
\end{lem}

\proof
%The proof essentially reduces to the calculations on page 308 of \cite{r5}.
By Lemma \ref{l3}
\beq\label{e1.12}
g_h^{k-1}(x)=\sum_{l=1}^k (-1)^{l-1} \binom{k-1}{l-1}
f(x+(l-1)h).
\e
Hence for $h=\frac{y-x}k$
\bml \label{e1.13}
\D_h^k f(x)=g_h^{k-1}(x)-g_h^{k-1}(x+h)\\
{}=\sum_{l=1}^k (-1)^{l-1}\binom{k-1}{l-1}[f(x+(l-1)h)-f(x+lh)]\\
{}=\sum_{j=0}^{k-1}(-1)^j\binom{k-1}j f(x+jh)
+\sum_{l=1}^{k}(-1)^{l}\binom{k-1}{l-1} f(x+lh)\\
{}=f(x)+(-1)^k f(x+kh)+\sum_{j=1}^{k-1}(-1)^j
\Bigl[\binom{k-1}j + \binom{k-1}{j-1}\Bigr] f(x+jh)\\
{}=(-1)^k\wt\D{}^k f(x,y)=\D^k f(x,y). \q{\square}
\e
%\endproof

Formulas \er{e4}--\er{e1.13} are examples of a series of formulas of finite
differences \cite{r8}, \cite{r7} which connect the operations of vector differential
operators $\n^lf$ with finite difference operators $\D_h^k f$ for various
values of the parameters $l,k$ and~$h$. They allow to reduce the estimates of
higher order difference remainders of \f s to lower order remainders of
higher order gradients of these \f s. They are analogues to the operations in
Taylor--Whitney's algebras which played an important role in the development
and applications of pointwise inequalities in \cite{Boj3}, \cite{Boj4}, 
\cite{r0}, \cite{r2} and in the Whitney--Glaeser--Malgrange theory
\cite{Gla1}, \cite{Gla}, \cite{Mal}, \cite{Wh1}, of smooth \f s on arbitrary
closed subsets of $\R^n$. Deep and important papers of Glaeser \cite{Gla2},
\cite{Gla3}, \cite{Gla} in this theory seem to be so far waiting for better
understanding and exploitment.

Now, combining Lemmata \ref{l1}--\ref{l4} we obtain the estimate \er{e1.4}, 
with $\wh a_f(x)$ controlled by the maximal \f\ of the vector 
gradient $\n^k f$ as required.

For convenience we summarize our discussion in the following

\begin{prop}\label{p1}
Let $f\in W^{m,p}(\R^n)$, $1<p\le\infty$. Then there exists a \f\ $\wh a_f\in
L_p(\R^n)$ such that the inequality \er{e1.4} holds for almost all points
$x,y\in \R^n$. The \f\ $\wh a_f$ is majorized a.e.\ by the local maximal \f\
of the generalized Sobolev gradient $|\n^m f|\in L_p(\R^n)|$.
\end{prop}

As formulated above Proposition \ref{p1} is the ``necessary part'' of the
main theorem in~\cite{r0} for the case of Lagrange interpolation remainders.
The proof of the ``sufficiency part'' of the theorem is left unchanged and
proceeds along the argument sketched in~\cite{r0}.

The sketch of the proof of Proposition \ref{p1} presented above is in the convention that
we work in the class of smooth \f s where all the operations involved are
classically meaningful. The clue of the story is that the constants appearing
in the estimates depend on the parameters $n$ and $p$ only. This fundamental
fact comes up again in~\SS2 below where the pointwise inequality appears as
``stable'' under smoothing by convolution.

\section{}
One of the advantages of the Lagrange \ip\ calculus and the Lagrange \r s
over the Taylor--Whitney \r s is that they interact very well with
convolutions. This is immediately seen for the \r s $\D^1f(x;y)=f(y)-f(x)$.

Indeed, for a normalized mollifier $\vf_\ve(x)$, $\vf_\ve\ge0$, $\int
\vf_\ve(x)\,dx=1$, we have
\[
f_\ve(x)\equiv f*\vf_\ve(x)=\int f(x-\eta)\vf_\ve(\eta)\,d\eta
\]
and
\[
f_\ve(y)-f_\ve(x)=\int \bigl[f(y-\eta)-f(x-\eta)\bigr]\vf_\ve(\eta)\,d\eta,
\]
hence
\[
|f_\ve(y)-f_\ve(x)|\le |x-y|
\Bigl(\int \wh a(x-\eta)\vf_\ve(\eta)\,d\eta +
\int \wh a(y-\eta)\vf_\ve(\eta)\,d\eta\Bigr)
\]
or
\beq \label{e2.1}
|f_\ve(y)-f_\ve(x)|\le |x-y|\bigl(\wh a_\ve(x)+\wh a_\ve(y)\bigr).
\e
This elementary though basic fact, already formulated in our paper \cite{r1}
and even earlier and repeated later in many seminar talks, holds for the higher order \r s
$\D^mf(x;y)$ as well. Indeed, e.g.\ for $m=2$ we have
\bmlg
f_\ve(x)-2f_\ve(\tfrac{x+y}2)+f_\ve(y)
=\int\bigl[f(x-\eta)-2f(\tfrac{x+y}2-\eta)+f(y-\eta)\bigr]\vf_\ve(\eta)\,d\eta\\
{}\le|x-y|^2 \int \bigl[\wh a_f(x-\eta)+\wh a_f(y-\eta)\bigr]\vf_\ve(\eta)
\,d\eta=|x-y|^2 \bigl(\wh a_{f,\ve}(x)+\wh a_{f,\ve}(y)\bigr)
\e
and the same calculation works for $m>2$.

Pointwise inequalities for the Lagrange \r s $\D^mf(x,y)$ appeared also in
the recent papers of H.~Triebel and his school \cite{HT}, \cite{Tri}. When
the paper \cite{r0} was written the papers of Triebel \cite{Tri} and
Haroske--Triebel \cite{HT} were unknown to the author. Geometrically for
$m>1$ they differ from ours by introducing the intermediate nodes $x_i=x+ih$,
$i=1,\dots,m-1$, in the right hand side of \er{e1.4}.

For arbitrary $0<p\le\infty$, $s>0$ and $m\in\mathbb N$ with $s\le m$ in \cite{Tri}
is introduced the class $\bL_p^{s,m}(\R^n)$ of all $f\in L_p(\R^n)$ for which
there exists a nonnegative \f\ $g\in L_p(\R^n)$ such that for all $h\in\R^n$,
$0<|h|\le1$ the inequality
\beq \label{e2.2}
|\D_n^m f(x)|\le |h|^s \sum_{l=0}^m g(x+lh) \q \hbox{a.e. in }\R^n
\e
holds.

With the norm
\[
\|f\|_{s,m,p}=\|f\|_{L_p(\R^n)}+\inf\|g\|_{L_p(\R^n)},
\]
where the infimum is taken over all $g$ admissible in \er{e2.2}, the space
$\bL_p^{s,m}(\R^n)$ is a quasi-Banach space \cite{Tri}.

For $s=m$, $p\ge1$, the inequality \er{e2.2} is stronger than \er{e1.4},
i.e.\ \er{e1.4} implies \er{e2.2}. Now it is clear that one of immediate
conclusions from the main theorem in~\cite{r0} is that for these values
of~$s$ and~$p$ \er{e1.4} is equivalent to~\er{e2.2}. Thus we conclude that the
intermediate nodes in~\er{e2.2} can be discarded.

For a while, let us introduce the notation: $\bW^{m,p}(\R^n)$ is the class of
all \f s $f\in L_p(\R^n)$ for which there exists an $\wh a_f\in L_p(\R^n)$
such that inequality \er{e1.4} holds a.e.\ in~$\R^n$.

In the general context of \er{e1.4} and \er{e2.2} we have the following

\begin{prop}\label{p6}
For $f\in \bW^{m,p}(\R^n)$ the mollified \f\ $f_\ve=f*\vf_\ve$ is in the class
$\bW^{m,p}(\R^n)$ and
\beq \label{e2.3}
|\D^m f_\ve(x;y)| \le |x-y|^m \bigl(\wh a{}_{f,\ve}^m(x)
+\wh a{}_{f,\ve}^m(y)\bigr)
\e
with $\wh a{}_{f,\ve}^m(x)=\wh a{}_f^m*\vf_\ve$.

Moreover, by the known properties of convolutions the $L_p$-norms of the mean
maximal gradients $\wh a_{f,\ve}^m$ are uniformly controlled by the
$L_p$-norms of mean maximal gradients of~$f$:
\beq\label{e2.4}
\|\wh a_{f,\ve}^m \|_{L_p(\R^n)}\le \|\wh a_{f}^m \|_{L_p(\R^n)}
\e
for all $\ve>0$.
\end{prop}

In particular, we conclude that smooth \f s are dense in $\bW^{m,p}(\R^n)$.

For the same values of the parameters $s,m,p$ as for $\bL_p^{s,m}(\R^n)$ we
can also consider the class $\wh\bL{}_p^{s,m}(\R^n)$ defined by the pointwise
inequality
\beq \label{e2.3/5}
|\D^mf(x;y)|\le |x-y|^s(g(x)+g(y)), \q y=x+mh,
\e
for some $g\in L_p(\R^n)$.

In \cite{Tri} the quasi-Banach spaces $\bL_p^{s,m}(\R^n)$ are used to
identify some Besov spaces $B_{p,q}^{s\theta}(\R^n)$, $0<\theta<1$, $0<q
\le\infty$, as real interpolation spaces (\cite{Tri}, the main theorem). For
$s=m$ in~\cite{r0}, as well as in the paper \cite{Tri}, the spaces
$\bW^{m,p}(\R^n)$ are identified with classical Sobolev spaces
$W^{m,p}(\R^n)$. 

The natural interesting question is to characterize the spaces $\wh
\bL{}_p^{s,m}(\R^n)$, for $s$ not integer, as some Sobolev--Besov type spaces. 

Proposition \ref{p6} and its proof are also valid for the class
$\wh\bL{}_p^{s,m}(\R^n)$. Thus smooth \f s are also dense in $\wh\bL_p^{s,m}
(\R^n)$.

The described characterization of Sobolev spaces $W^{m,p}(\R^n)$, $p>1$, by
pointwise inequalities \er{e1.4} obviously holds for open subdomains
$G\subset \R^n$ as well, if they have sufficiently regular boundary, e.g.\ for
extension domains \cite{Shv2}. However this characterization is definitely
not true for arbitrary subdomains. In this context it seems legitimate to
introduce the (maximal) class of subdomains which admit global pointwise
characterization. 

\begin{dfn}
An open subdomain $G\subset\R^n$ is called a \ti{natural Sobolev $(p,s)$,
$1<p\le\infty$, $s>0$, domain} if the pointwise inequality
\beq \label{e2.5}
|\D^mf(x,y)|\le |x-y|^s (a(x)+a(y)), \q s\le m,
\e
for all pairs of points $x,y\in G$ such that the segment $[x,y]\subset G$,
defines a Sobolev type Banach space $\wt W{}^{s,p}(G)$.
\end{dfn}

For $s\le m$, $0<p\le \infty$ a related class of spaces $L_p^s(\R^n)^m$ has been introduced by
H.~Triebel in~\cite{Tri} and identified with subspaces of Besov spaces
$B_{p,\infty}^s (\R^n)$ for $G=\R^n$.

It is natural to ask in particular in what sense our Sobolev--Besov type
spaces coincide with Triebel--Besov spaces in~\cite{Tri}. More generally we
can ask in what sense and for what subdomains the spaces $\wt W{}^{s,p}(G)$
coincide with Besov--Sobolev type spaces for non-integer~$s$.

\section{}
The pointwise inequality \er{e1.4} actually can be used to characterize
rather the homogeneous Sobolev spaces $\dot W{}^{m,p}(G)$ for a subdomain
$G\subset \R^n$ with the seminorm $\|\nabla^m f\|_{L^p(G)}$. The classical
inhomogeneous Sobolev spaces arise then as subspaces of $L^p(G)$. For the model
case $G=\R^n$ various delicate phenomena are related to the asymptotic
behavior of \f s in $\dot W{}^{m,p}(\R^n)$ for $|x|\to\infty$ which seem to
be so far only partially understood (see \cite{BIN} and numerous other
references. See also \cite{Ba} for the case $n=1$, i.e.\ on the line $\R^1$).

In our presentation here the pointwise inequality \er{e1.4} as well as
somewhat more sophisticated inequality for Taylor--Whitney \r s
$R^{m-1}f(x;y)$ should be considered as elementary, though fundamental, facts
appearing at the first steps of any discussion of Sobolev spaces.

What seems to be still lacking in this elementary discussion of Sobolev space
theory, is the deeper, geometric and analytical, understanding of the trace
(projection) operator $W^{m,p}(\R^n)\to W^{s,p}(\R^k)$ for the corresponding
values of the parameters $s,p$. As is classically known, this question led to
the introduction of Sobolev fractional spaces, $W^{s,p}(\R^n)$, $s$~---~real,
$\cH^{s,p}(\R^n)$ and Besov spaces (\cite{BIN}, \cite{St1}, \cite{St2} and
many references therein).

It seems also that special attention should be directed to ``\ti{les
sch\'emas d'interpolation}'' of G.~Glaeser \cite{Gla1}, \cite{Gla},
\cite{Ker}, \cite{MM}, and their role in the theory of Sobolev spaces, see also
\cite{Shv1}, \cite{Shv2}.

As already remarked in \cite{r0}, and earlier even in \cite{Boj3}, \cite{r1},
\cite{r2}, the pointwise inequality \er{e1.4} above, together with the
related inequality for the Taylor--Whitney \r\ $R^{m-1}f(x;y)$ (precisely,
inequality (1.2) in~\cite{r0}) may serve as natural starting points and
effective tools in the study of fundamental structural properties of Sobolev
\f s. Let us briefly recall some of them without going here into details
(postponed to the activity and exposition plan foreseen in the last lines
of~\cite{r0}). 

\def\labelenumi{\alph{enumi})}
\begin{enumerate}
\item Lusin's approximation of Sobolev \f s, i.e.\ interpolation by smooth \f
s on closed
subsets, up to open complements of arbitrary small measure.
\item Stability under the convolution with compactly supported $C^\infty$
kernels; density of subspaces of smooth, $C^\infty$ \f s.
\item S. M. Nikolskii's \cite{BIN} fundamental theorems describing the
characterization of Sobolev \f s by their behavior on typical (almost all, in
some natural sense) hyperplanes of positive codimensions less than $n$. The
inequalities \er{e1.4} reduce the characterization of Sobolev \f s to their
behavior on affine segments in their domain of definition. In particular the
\f s in $W^{m,p}(\R^n)$ for $p>1$ are H\"older continuous on a.e.\
hyperplanes $\R^k\subset \R^n$, $k<p$, the only global condition binding the
variable Lipschitz \cf s on hyperplanes is the Fubini theorem for Lebesgue
spaces $L^p(\R^n)$ and the factorizations $\R^n \sim \R^{n-k}\times \R^k$.
\item Characterization of compact subsets of $W^{m,p}(\R^n)$ by some
conditions on the corresponding variable Lipschitz \cf s.
\item Differentiability properties of Sobolev \f s, Calder\'on
differentiability theorems, approximate and Peano differentiability,
\cite{BIN}, \cite{r1}, \cite{r2}, \cite{DL}, \cite{DS}, \cite{Ha2}, \cite{MZ}, 
\cite{Rad}, \cite{Stp}.
\item Extension of Haj\l asz--Sobolev imbedding theorems of Sobolev spaces
into higher exponent Lebesgue spaces: $W^{l,p}(G)\subset L^q(G)$, $q>p$,
$l>1$, for suitable values of the parameters $l,p,q$ and $\dim G=n$, as in
the classical Sobolev theory, modeled on Haj\l asz's proof for measure metric
spaces, \cite{Ha1}.
\item Extension of the classical Hermite \ip\ formulas \cite{r8}, \cite{r7}
to the general multidimensional context of ``multiple'' nodes as a
theory intermediate between the Taylor--Whitney and Lagrange (``simple''
nodes) \ip\ theory (\cite{Gla1}, \cite{Gla}, \cite{Mal}).
\item If instead of the colinear equidistant nodes other configurations of
\ip\ nodes are used, more complicated algebraic and geometric phenomena
occur, e.g.\ in the case of colinear not equidistant simple nodes the
classical divided difference calculus comes up \cite{DS}, \cite{Fi}. New \ip\
concepts appear in the Glaeser papers \cite{Gla1}, \cite{Gla3}, and, more
recent, \cite{Ha1}, \cite{Mal}. Combining these ideas with the Sobolev's
averaging procedure, \cite{Sob}, \cite{Sob3}, \cite{Sob2}, 
defines an apparently new interesting research direction.
\item Also the natural inclusions, e.g.\ $W^{m+1,p}(\R^n)\subset
W^{m,p}(\R^n)$ when interpreted in terms of pointwise inequalities \er{e1.4}
suggest direct implications between inequalities \er{e1.4} for various
admissible values of the parameters $(m,p)$. These lead to interesting and
non-trivial arguments of geometric and analytic character. Probably the first
beautiful example of this type of argument was given by Y.~Zhou \cite{Zh}.
\end{enumerate}

Let us remark at last that the ``pointwise'' approach to the Sobolev space
theory seems to bring out more clearly than usually presented in the
literature (e.g.\ \cite{Kra}, \cite{Mar}, \cite{MZ}, \cite{r8}, \cite{DL})
close and natural connections between the general concepts of \f al spaces
and approximation theory on the real line $\R^1$ and in $n$-dimensional,
$n>1$, euclidean spaces.

\end{document}